\documentclass[letterpaper,times]{sae}

\usepackage[svgnames]{xcolor}
\definecolor{SAEblue}{RGB}{1,160,233}
\raggedright

\usepackage[justification=justified, singlelinecheck=false]{caption}  
\DeclareCaptionFont{eightpt}{\fontsize{8pt}{11pt}\selectfont #1}
\usepackage{lineno}
\usepackage[utf8]{inputenc}

\usepackage{array}
\newcolumntype{L}[1]{>{\raggedright\let\newline\\\arraybackslash\hspace{0pt}}p{#1}}
\newcolumntype{C}[1]{>{\centering\let\newline\\\arraybackslash\hspace{0pt}}p{#1}}
\newcolumntype{R}[1]{>{\raggedleft\let\newline\\\arraybackslash\hspace{0pt}}p{#1}}

\usepackage{graphicx}
\usepackage{multirow}
\usepackage{amsmath}
\usepackage{amssymb}
\usepackage{comment}

\newcommand{\ignore}[1]{}

\usepackage{nameref}
\usepackage{booktabs}

\makeatletter
\def\@seccntformat#1{%
  \expandafter\csname c@#1\endcsname\c@section
  }
\makeatother

\usepackage{subfigure}
\usepackage{multimedia}
\usepackage{lipsum}

\PaperTitle{A Decentralized Time- and Energy-Optimal Control Framework for Connected Automated Vehicles: From Simulation to Field Test}
\usepackage{calc} 

\makeatletter 
\renewcommand\@biblabel[1]{#1. } 
\makeatother

\AddAuthor{A M Ishtiaque Mahbub}{University of Delaware, USA}
\AddAuthor{Vasanthi Karri, Darshil Parikh, Shyam Jade}{Robert Bosch LLC., USA}
\AddAuthor{Andreas A. Malikopoulos}{University of Delaware, USA}

\PaperNumber{20XX-XX-XXXX} 
\SAECopyright{2020} 

\begin{document}
\maketitle
\section{Abstract}
The implementation of connected and automated vehicle (CAV) technologies enables a novel computational framework for real-time control aimed at optimizing energy consumption with associated benefits. In this paper, we implement an optimal control framework, developed previously, in an Audi A3 etron plug-in hybrid electric vehicle, and demonstrate that we can improve the vehicle's efficiency and travel time in a corridor including an on-ramp merging, a speed reduction zone, and a roundabout. Our exposition includes the development, integration, implementation and validation of the proposed framework in (1) simulation, (2) hardware-in-the-loop (HIL) testing, (3) connectivity enabled virtual reality based bench-test, and (4) field test in Mcity. We show that by adopting such inexpensive, yet effective process, we can efficiently integrate and test the control framework, establish proper connectivity and data transmission between different modules of the system, and reduce uncertainty. We evaluate the performance and effectiveness of the control framework and observe significant improvement in terms of energy and travel time compared to the baseline scenario.

\section{Introduction}

The increasing traffic volume in urban areas has reached the capacity of current infrastructure resulting in congestion. Fuel efficiency and travel time can also be seriously affected in daily commute \cite{Schrank2015}. 
Connected and automated vehicles (CAVs) provide the most intriguing opportunities for enabling users to better monitor transportation network conditions and make better operating decisions. The interconnectivity of mobility systems enables a novel computational framework to process massive amount of data, and deliver real-time control actions that optimize energy consumption and associated benefits. The advent of machine-to-machine (M2M) or internet of things (IoT) connectivity protocol and their use in vehicular systems has paved the way for emerging technologies like vehicle-to-everything (V2X) communication systems, which is the first step to autonomous driving and connected road infrastructure. The IoT, which is a system of interrelated computing devices enabling data transfer over a network without requiring human interaction, has transformed the CAVs into moving data centers. From a control and optimization view-point, CAVs can alleviate congestion at the major transportation segments such as urban intersections, merging roadways, roundabouts, and speed reduction zones, which are the primary sources of bottlenecks that contribute to traffic congestion. Several research efforts can be found in the literature proposing centralized, or decentralized control approaches to alleviate congestion at these bottlenecks in order to improve fuel efficiency and travel time. A detailed discussion of these research efforts reported in the literature to date can be found in \cite{Malikopoulos2016a}. In earlier work, a decentralized optimal control framework was established for coordinating online CAVs in different transportation segments. A closed-form, analytical solution was presented in \cite{Rios-Torres2} and \cite{Ntousakis:2016aa} for coordinating online CAVs at highway on-ramps, in \cite{Malikopoulos2017} and in \cite{Mahbub2019ACC} at intersections, in \cite{Malikopoulos2018a} at roundabouts, in \cite{Malikopoulos2018c} for speed harmonization, and in \cite{Zhao2018} in a corridor with several conflict zones. The solution of the optimal control problem considering state and control constraints was presented in \cite{Malikopoulos2017} and \cite{Mahbub2019CDC} at an urban intersection without considering rear-end collision avoidance constraint, and the conditions under which the latter does not become active were presented in \cite{Malikopoulos2018c}. These efforts considered a control zone inside of which the CAVs can communicate with each other through vehicle-to-vehicle (V2V) connectivity, and with a coordinator through vehicle-to-infrastructure (V2I) connectivity, giving rise to the decentralized communication topology. The coordinator is not involved in any decision making process. It only handles the communication of appropriate information among CAVs. The performance and the effectiveness of the aforementioned decentralized vehicle dynamics (VD) controllers for individual tasks (i.e., on-ramp merging, roundabout, signal-free intersection) have been validated at University of Delaware's 1:25 Scaled Smart City  \cite{beaver2019demonstration} using robotic cars. The necessity to validate the effectiveness of the proposed VD controllers in real-world test environment using actual vehicles is the next natural step.

Currently, CAV testing and performance evaluation can be categorized into the following conventional approaches: testing (1) in a simulation environment, (2) at a closed and controlled test facility, and (3) on public roads. These approaches have associated advantages and drawbacks. 
Vehicle testing on public roads has obvious limitations \cite{Singhvi2016}. Moreover, it is often important to conduct system validation to find the \textit{corner cases} that can challenge the autonomous vehicles. The corner cases are rare events often associated with the exception handling of the controller algorithms or vehicle operation at the extremes of acceptable operating parameters. In the real world traffic, these corner cases are not easy to capture due to the presence of many external variables. As a consequence, although several companies have been demonstrating their self-driving cars on public roads, the topic of whether or not to allow CAVs to run in conjunction with general traffic raises significant debate \cite{Fagnant2015}. Vehicle testing in a closed-test facility is also associated with safety issues, wastage of fuel, and high resource and maintenance costs. On the contrary, simulation is a cost-effective way to test new technologies, and is safer, cheaper and faster compared to using test vehicles. Simulators can be very important tools, as they can help us gather key information about consumer's preferences and behaviors. However, testing driverless vehicles in a simulation environment may not be as accurate as testing in real world due to the difficulty in modeling the exact physics of the test environment, vehicle dynamics and driving behavior in simulation. The complexities arising from data loss, transmission latency, etc, associated with the V2X communication framework cannot be captured either. Therefore, an intermediate measure has to be employed which can subscribe to the benefits of both the simulation environment as well as the real-world testing procedure. Some studies developed hardware-in-the-loop (HIL) \cite{Zulkefli2017} simulation platform, where an engine with virtual powertrain model was used. Another approach developed vehicle-in-the-loop \cite{Sunkari2016}  simulation platform which used an entire vehicle into the simulation. To model real vehicle behavior observed in the field, a parallel traffic system was proposed, which set up a mirror of the real world in virtual space \cite{Wang2010_its, Li2016a}. The parallel system can be used to design different test scenarios and evaluate how vehicles perform in these scenarios \cite{Li2016b}. More recently, test procedure based on virtual reality has been successfully implemented in Mcity vehicle testing facility. In this framework, the performance of the vehicle is tested in the real world facility, and additional traffic congestion in the network is created using virtual vehicles from a simulator. The real-world test vehicle can interact with the virtual vehicles through V2X communication, and the behavior of the virtual vehicles and the test vehicle are synchronized \cite{yiheng2019}.

To harness the advantages of the aforementioned vehicle testing methods, we adopt a systematic and sequential approach to validate the proposed VD controller in a fixed traffic corridor including an on-ramp merging, a speed reduction zone, and a roundabout. We follow a sequential approach, e.g., from simulation environment to the real-world field test, to validate the effectiveness of the VD controller in the corridor. First, we investigate the effectiveness of the controller in a commercial simulation environment, and evaluate its performance in terms of energy consumption and travel time reduction. Then, we make necessary modifications (hardware/software) to our test vehicle, which is an Audi A3 e-tron plugin hybrid electric vehicle (PHEV), and re-evaluate the controller performance in a HIL test through a chassis-dyno setup. In this context, we integrate the test vehicle with the proposed controller and associated communication protocol, and conduct system validation by setting up a bench test model. Finally, we run a virtual reality enabled field test in Mcity with full V2X communication capability, and quantify the real-world implication of the proposed vehicle dynamics controller.


%
In this paper, we (1) implement a vehicle coordination control algorithm for throughput maximization and energy optimization in the corridor without any traffic signals, (2) validate and evaluate the performance indices through a simulated environment and hardware-in-the-loop test facility, (3) present an IoT-enabled virtual reality based bench test model to perform system validation, and (4) evaluate the performance of the VD controller in a V2X enabled virtual reality based real-world test facility in Mcity.\\
%

The remainder of the paper is organized as follows: First, we present the problem formulation and modeling framework, and derive the analytical, closed form solution of the VD controller. Then, we validate the effectiveness of the efficiency of the VD controller, and present a sequential approach including a simulation environment, a HIL test, and a bench test. Finally, we provide a detailed exposition into the virtual reality based field test with full V2X functionality in Mcity, discuss the results, and draw concluding remarks.

\section{Problem Formulation}
\begin{figure}[h]
\centering
\includegraphics[scale=0.75]{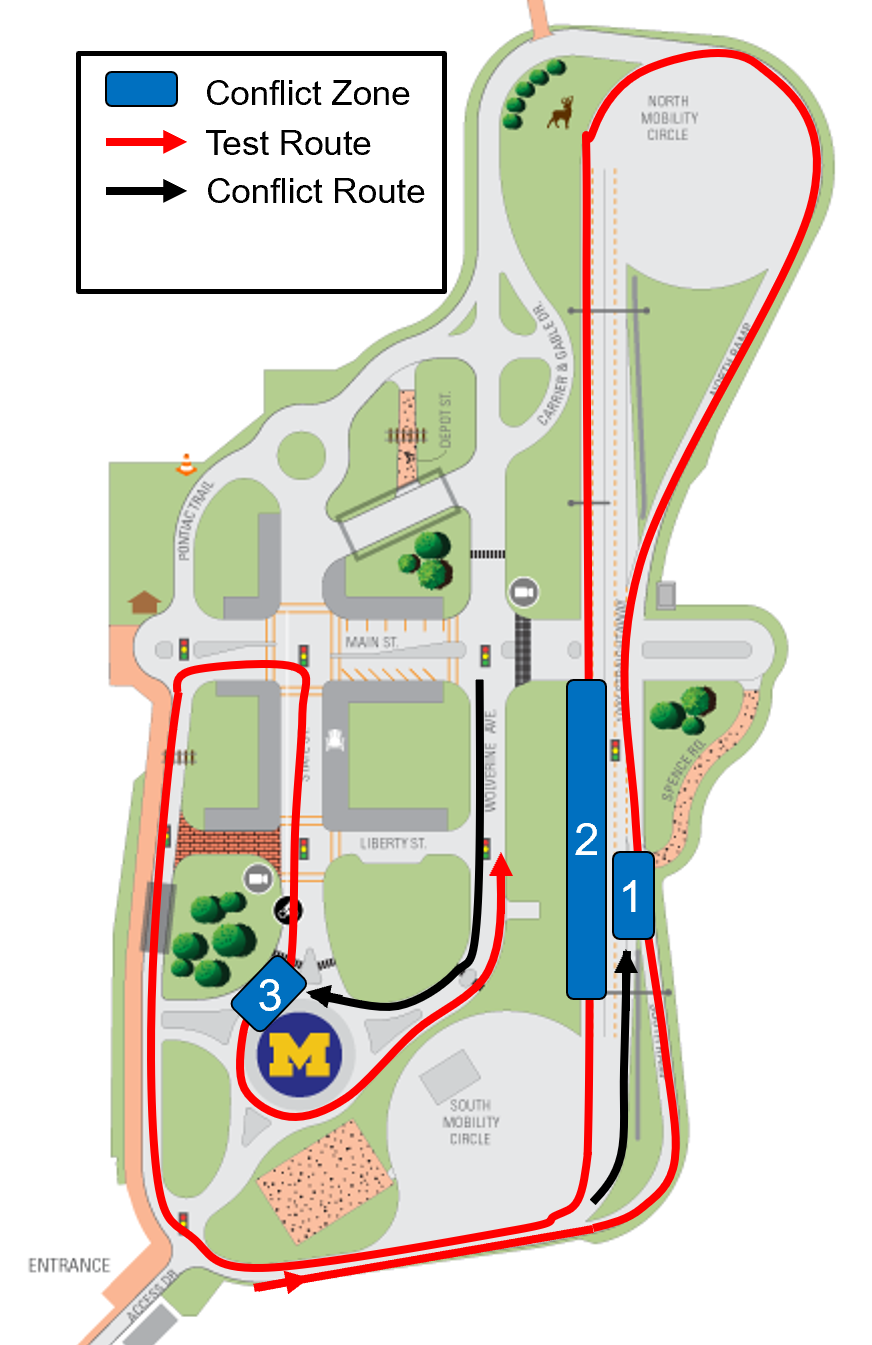} 
\caption{The corridor in Mcity with three conflict scenarios: (1) on-ramp merging, (2) speed reduction zone, and (3) roundabout.}%
\label{fig:mcity_corridor}%
\end{figure}

We consider a network of CAVs driving through a corridor of Mcity, which is a 32-acre vehicle testing facility located at the North campus of University of Michigan. A typical commute of a vehicle usually includes merging at roadways, crossing signalized intersections, cruising in congested traffic, passing through speed reduction zones, etc. Therefore, the VD controller needs to be able to optimize the speed profile of the vehicle over such traffic scenarios. To validate the performance of the VD controller, we consider a particular corridor of the Mcity traffic network that consists of several conflict scenarios, e.g., a highway on-ramp, a speed reduction zone (SRZ), and a roundabout as shown by the numbered blue rectangles in Fig. \ref{fig:mcity_corridor}. Each CAV follows the test route illustrated by the red trajectory in Fig.  \ref{fig:mcity_corridor} through the conflict scenarios $z=1,2,3$ denoting respectively the highway on-ramp, the SRZ, and the roundabout. Note that, to create traffic congestion in the test route, we consider additional \textit{conflict routes} as shown by the black trajectories in Fig. \ref{fig:mcity_corridor}. 

\begin{figure}[h]
\centering
\includegraphics[width=3.5in]{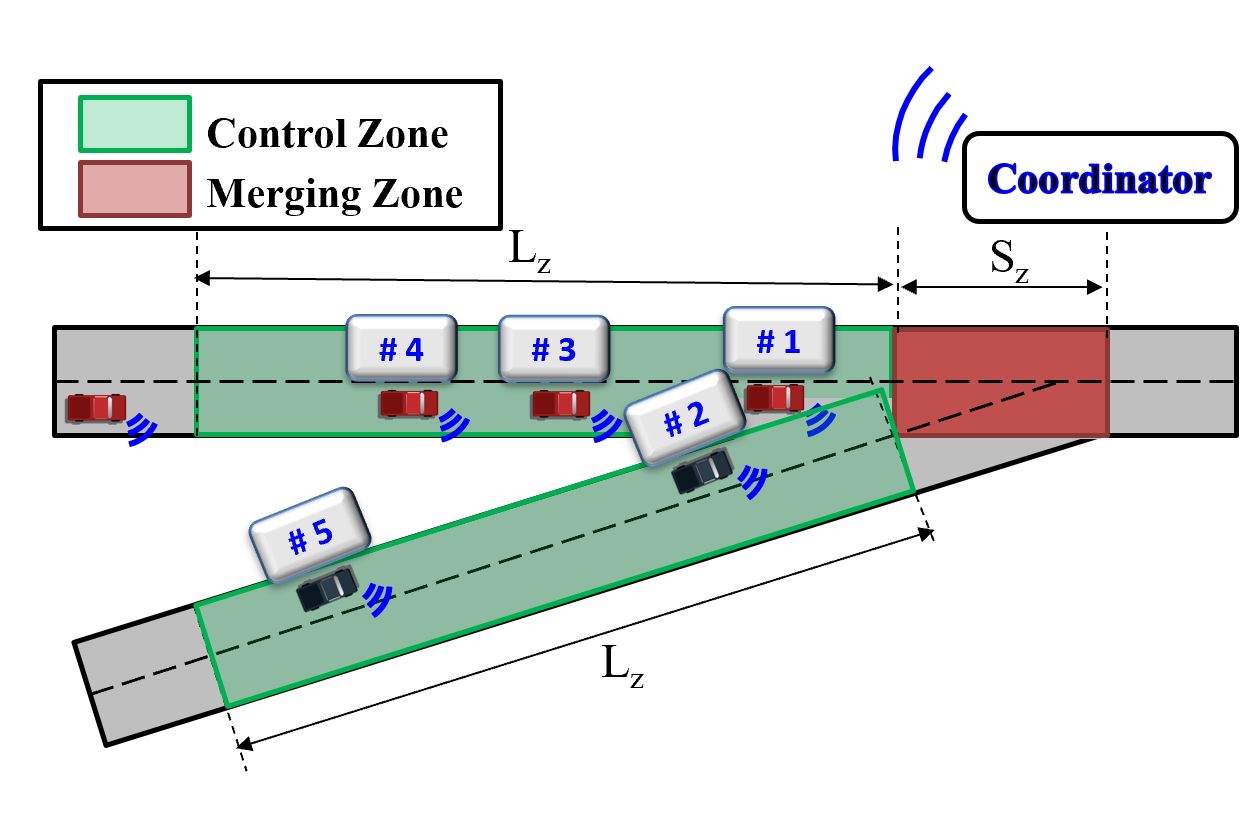} 
\caption{An automated on-ramp merging scenario with coordinator.}%
\label{fig:ramp_merging}%
\end{figure}

In Fig. \ref{fig:ramp_merging}, an automated on-ramp merging scenario with CAVs is illustrated. The essential features of this conflict scenario can be generalized to the other conflict scenarios, i.e., the SRZ and the roundabout. Upstream of each conflict scenario, we define a \textit{control zone} (CZ) where CAVs coordinate with each other to avoid any rear-end or lateral collision in the conflict zone. The length of the CZ $L_z$ specific to a conflict scenario $z$ can be varied, and is restricted by the range of the $coordinator$ present (see Fig. \ref{fig:ramp_merging}). For each conflict scenario, the coordinator communicates with the CAVs traveling within its control range. Note that, the coordinators do not make any control decisions for the CAVs. They only handle the information among the CAVs and maintain the queue that each CAV enters the corresponding CZ. For example, in Fig. \ref{fig:ramp_merging} we observe that the CAVs inside the CZ have been assigned a unique identity, whereas the vehicle outside the CZ is yet to be communicated by the coordinator. We define the area of potential lateral collision to be the $merging~zone$ (MZ) as illustrated by the red square in Fig. \ref{fig:ramp_merging}. We denote the length of the MZ to be $S_z$ specific to conflict scenario $z$.
The objective of each CAV is to derive its optimal control input (acceleration/deceleration) to cross the conflict zones without any rear-end or lateral collision with the other vehicles. The VD controller generates optimal speed trajectory to simultaneously minimize the travel time and energy consumption over different traffic scenario. 

In the modeling framework described above, we impose the following assumptions:

\textbf{Assumption 1: }For each CAV $i$, none of the constraints is active while entering a CZ of conflict scenario $z,~z= 1,2,3$.

\textbf{Assumption 2: }Each vehicle is equipped with sensors and V2X communication device to measure and share their local state information. 

\textbf{Assumption 3: }The communication among CAVs occurs without any transmission latency, errors or data loss. 

\textbf{Assumption 4: }No lane change is allowed inside the vehicle route.

\textbf{Assumption 5: }The speed of each CAV $i$ inside the merging zone is constant.

The first assumption ensures that the initial state and control input of the CAVs while entering the CZ are feasible. The second assumption enables mathematical formulation of the VD controller with $100\%$ CAV penetration rate. The third assumption might be strong, but can be relaxed as long as the noise in the measurements and/or delays is bounded. For example, we can determine upper bounds on the state uncertainties as a result of sensing or communication errors and delays, and incorporate these into more conservative safety constraints. The fourth assumption simplifies the formulation by restricting the traffic flow to a single lane. Hence, only longitudinal dynamics of the vehicle is considered. Finally, the fifth assumption enables constant space headway between two consecutive CAVs to avoid lateral and rear-end collision within the MZ.

\begin{figure}[h]
\centering
\includegraphics[width=3.5in]{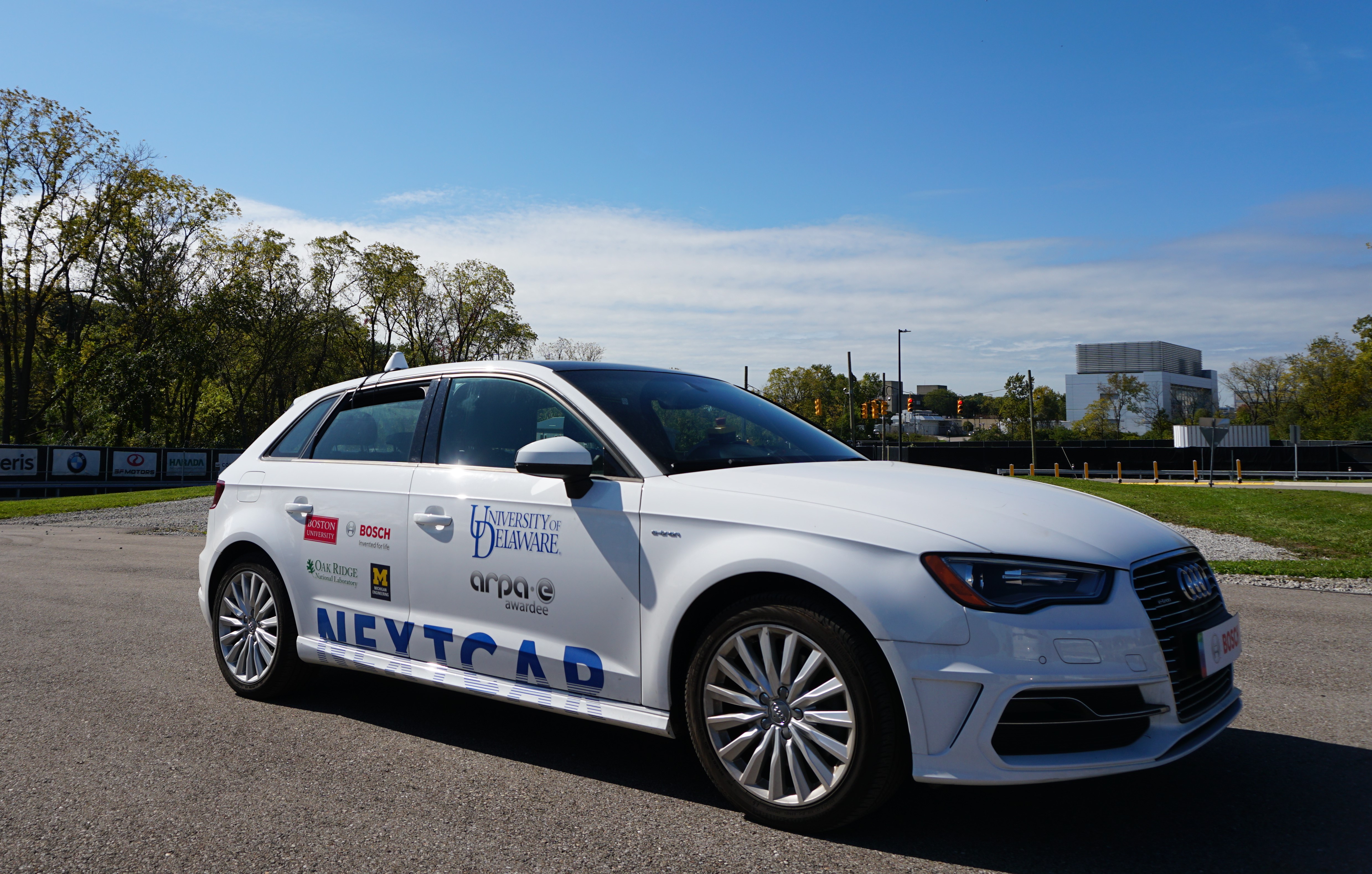} 
\caption{The test vehicle, Audi A3 e-tron, in Mcity.}%
\label{fig:audi}%
\end{figure}

Our goal is to validate the effectiveness of the VD controller and investigate its performance in a real-world scenario. To this end, we adopt a methodology which includes sequential steps, and systematically proceed from one to another. We start with the controller validation in the simulation environment, and end up conducting the real-world field test in Mcity with the help of virtual reality. These sequential steps can be briefly introduced as following:
\begin{enumerate}
    \item \textbf{Simulation Environment:} We create a simulation environment of Mcity in a commercial software (i.e., PTV VISSIM), and develop the baseline and optimal scenarios.The baseline scenario considers the case of 0\% penetration rate of CAVs. The optimal scenario corresponds to the 100\% penetration rate of CAVs, and they all are controlled by the VD controller.
    \item \textbf{Hardware-in-the-Loop Test:}  We make necessary modifications to the Audi A3 in terms of vehicle hardware and software to be able to implement the VD controller output (recommended speed).  We feed the baseline and optimal controlled speed profiles generated in the simulation environment in step 1 to the Audi A3, and trace them in a chassis-dyno setup. 
    \item \textbf{Bench Test:} We integrate the controller in the Audi A3's head unit. We develop a bench-test method and conduct a basic system validation. Here, we employ the message queuing telemetry transport (MQTT) protocol, which is a publish-subscribe based messaging protocol of IoT, in conjunction to a simulation based virtual traffic environment. We establish both the controller integration and the V2X communication framework.
    \item  \textbf{Field Test:} After following these steps, the Audi A3 integrated with the proposed VD controller can be taken to the test facility. Steps 1-4 not only reduce the possibility of technical complications during the field test but also provide an extra layer of safety to the driver and the test vehicle components. Therefore, these prerequisite steps can be considered essential before the on-field vehicle test session. Finally, we conduct the field test in a V2X enabled virtual reality based test facility in Mcity with the Audi A3 (Fig. \ref{fig:audi}). We provide a detailed exposition into the various aspects of the field test in the following sections.
\end{enumerate}

\section{Optimal Vehicle Dynamics Controller}
\subsection{Vehicle Dynamics and Constraints}
In our control approach, we consider a corridor that contains three conflict scenarios (Fig. \ref{fig:mcity_corridor}), e.g., a merging roadway, a SRZ, and a roundabout. Let $z= 1,2,3$ be the index of a conflict scenario in the corridor.  When a vehicle enters the CZ of a specific conflict scenario, the coordinator receives its information and assigns a unique identity $i$ to the vehicle. Let $t_i^{0,z}$ be the time when vehicle $i$ enters the CZ towards conflict scenario $z$, $t_i^{m,z}$ be the time when the vehicle enters the MZ of the conflict scenario $z$, and $t_i^{f,z}$ be the time when vehicle $i$ exits the corresponding CZ. Let $\mathcal{N}_z(t)=\{1,\dots,N(t)\}$, where $ t\in \mathbb{R}^+$ is the time, be a queue of vehicles to be analyzed for conflict scenario $z$. Here, $N(t)$ is the total number of CAVs within the CZ of the specific conflict scenario $z$ at time $t \in \mathbb{R}^+$. The dynamics of each vehicle $i, i\in \mathcal{N}_z(t)$, are represented with a state equation
\begin{equation} \label{eq:state}
\dot{x}(t) = f(t, x_i, u_i), ~ x_i(t_i^{0,z}) = x_i^{0,z},
\end{equation}
where $x_i(t), u_i(t)$ are resepectively the state and control input of the vehicle. We define $x_i(t_i^{0,z})=x_i^{0,z}$ to be the value of the initial state corresponding to the initial time $t_i^{0,z}$ when the vehicle $i$ enters the CZ of conflict scenario $z$. For simplicity, we model each vehicle as a double integrator, i.e., $\dot{p}_i = v_i(t)$ and $\dot{v}_i = u_i(t)$, where $p_i(t) \in \mathcal{P}_i, v_i(t) \in \mathcal{V}_i$, and $u_i(t) \in \mathcal{U}_i$ denote the position, speed, and acceleration/deceleration (control input) of each vehicle $i$. Let $x_i(t)=[p_i(t)~ v_i(t)]^T$ denote the state of each vehicle $i$, with initial value $x_i(t_i^{0,z})=[p_i(t_i^{0,z}) ~ v_i(t_i^{0,z})]^T$, taking values in the state space $\mathcal{X}_i=\mathcal{P}_i\times  \mathcal{V}_i$. The sets $\mathcal{P}_i, \mathcal{V}_i$, and $\mathcal{U}_i, i\in \mathcal{N}_z(t)$, are complete and totally bounded subsets of $\mathbb{R}$. The state space $\mathcal{X}_i$ for each vehicle $i$ is closed with respect to the induced topology on $\mathcal{P}_i\times  \mathcal{V}_i$ and thus, it is compact.

To ensure that the control input and vehicle speed are within a given admissible range, the following constraints are imposed.
\begin{equation}\label{eq:constraints}
\begin{aligned} 
u_{min} &\leq u_i(t) \leq u_{max},~ \text{and} \\
0 &\leq v_{min} \leq v_i(t) \leq v_{max}, ~ \forall t \in [t_i^{0,z}, ~ t_i^{f,z}],
\end{aligned} 
\end{equation}
where $u_{min}, u_{max}$ are the minimum deceleration and maximum acceleration respectively, $v_{min}, v_{max}$ are the minimum and maximum speed limits respectively. We assume homogeneity in terms of vehicle types, which enables the use of constant maximum and minimum acceleration values for any vehicle $i$.

To avoid the rear-end collision of two consecutive vehicles traveling on the same lane, the position of the preceding vehicle should be greater than, or equal to the position of the following vehicle plus a predefined safe distance $\delta_i(t)$, where $\delta_i(t)$ is proportional to the speed of vehicle $i$, $v_i(t)$. Thus, we impose the rear-end safety constraint
\begin{equation} \label{eq:safety}
s_i(t)=p_k(t)-p_i(t)\geq \delta_i(t), \forall t \in  [t_i^{0,z}, ~ t_i^{f,z}],
\end{equation}
where vehicle $k$ is immediately ahead of $i$ on the same lane. 

For each CAV $i\in \mathcal{N}_z(t)$, the lateral collision is possible within the set $\Gamma_i$,
\begin{equation}
\Gamma_i \overset{\Delta}{=} \{t\,\,|\,t\in  [t_i^{m,z}, t_i^{f,z}]\}.
\label{eq:lateral_active}
\end{equation}
Lateral collision between any two CAVs $i,j\in \mathcal{N}_z(t)$ can be avoided if the following constraint hold,
\begin{equation}
\Gamma_i \cap \Gamma_j=\varnothing, ~  \forall t\in [t_i^{m,z}, t_i^{f,z}], \quad i,j\in \mathcal{N}_z(t).
\label{eq:lateral_constraint}
\end{equation}

\subsection{Hierarchical Control Structure}
We formulate the following optimal control problem for each vehicle in the queue $\mathcal{N}_z(t)$
\begin{gather} \label{eq:min}
\min_{u_i}\frac{1}{2}\int_{t_i^{0,z}}^{t_i^{m,z}} u_i^2(t)dt, ~\forall i \in \mathcal{N}_z(t), ~\forall z =1,2,3, \\
\text{subject to}: (\ref{eq:state}), (\ref{eq:constraints}),\nonumber\\
p_{i}(t_i^{0,z})=p_{i}^{0,z}\text{, }v_{i}(t_i^{0,z})=v_{i}^{0,z}\text{, }p_{i}(t_i^{m,z})=p_{z},\nonumber\\
\text{and given }t_i^{0,z}\text{, }t_i^{m,z},\nonumber
\end{gather}
where $p_{z}$ is the location (i.e., entry position) of MZ $z$, $p_{i}^{0,z}$, $v_{i}^{0,z}$ are the initial position and speed of vehicle $i$ when it enters the CZ of conflict scenario $z$. In Eq. \eqref{eq:min}, we consider the cost function to be the $L^2$-norm of the control input, i.e., $\frac{1}{2}u_i^2(t)$, to minimize transient engine operation. By minimizing transient engine operation, we have direct benefits in fuel consumption. The solution of this problem formulation yields a closed-form analytical solution that can be implemented on-board the CAV in real-time.
To address this problem, we apply the decentralized optimal control framework and analytical solution presented in \cite{Malikopoulos2017}. Note that, in \eqref{eq:min} we do not explicitly include the lateral \eqref{eq:lateral_constraint} and rear-end \eqref{eq:safety} safety constraints. We enforce the lateral collision constraint by selecting the appropriate merging time $t_i^{m,z}$ for each CAV $i$ by solving the upper-level vehicle coordination problem. The activation of rear-end safety constraint can be avoided under proper initial conditions $[t_i^{0,z},v_i(t_i^{0,z})]$ as discussed in \cite{Malikopoulos2018c}.  

We define two sets $\mathcal{L}_i^z(t)$ and $\mathcal{C}_i^z(t)$ which contains the preceding vehicle travelling in the same lane, or in the conflict lane relative to CAV $i$ respectively. For example, CAV $\#3\in\mathcal{L}_4^1(t)$, and CAV $\#1\in\mathcal{C}_2^1(t)$ according to Fig. \ref{fig:ramp_merging}. Each vehicle $i$ determines the time $t_i^{m,z}$ that will be entering the conflict zone $z= 1,2,3,$ upon arrival at the entry of the corridor. The next vehicle $(i+1)$, upon its arrival at the entry of the corridor, will search for feasible times to cross the conflict zones based on available time slots. Therefore, each CAV $i$ recusively uses the information of the previous vehicle $(i-1)$, and calculates its time to enter the conflict zone.

If CAV $(i-1)\in \mathcal{L}^z_i(t)$,
\begin{equation}
    t_i^{m,z} = \max\left \{\min\left\{t_{i-1}^{m,z}+\frac{\delta(v_i(t))}{v_{i-1}(t_{i-1}^{m,z})}, \frac{L_z}{v_{\min}}\right \}, \frac{L_z}{v_0(t_i^{0,z})}, \frac{L_z}{v_{\max}}\right \}. \label{eq:7}
\end{equation}
If CAV $(i-1)\in \mathcal{C}^z_i(t)$,
\begin{equation}
    t_i^{m,z} = \max\left \{\min\left\{t_{i-1}^{m,z}+\frac{S_z}{v_{i-1}(t_{i-1}^{m,z})}, \frac{L_z}{v_{\min}}\right \}, \frac{L_z}{v_0(t_i^{0,z})}, \frac{L_z}{v_{\max}}\right \}. \label{eq:8}
\end{equation}
The recursion is initialized when the first vehicle enters the CZ, i.e., it is assigned $i = 1$. The details of upper-level optimization and the closed-form solution of the low-level optimal control problem for every vehicle, i.e., the optimal acceleration profile that will achieve the pre-determined  entry time at each conflict scenario, was presented in \cite{Zhao2018}.

%

To derive the analytical solution of \eqref{eq:min} using Hamiltonian analysis, we follow the standard methodology used in optimal control problems with control and state constraints \cite{bryson1975applied}. 
The solution of the constrained problem has been addressed in \cite{Malikopoulos2017}. If the unconstrained solution violates any of the state or control constraints, then the unconstrained arc is pieced together with the arc corresponding to the violated constraint, and we re-solve the problem with the two arcs pieced together. The two arcs yield a set of algebraic equations which are solved simultaneously using the boundary conditions of and interior conditions between the arcs. If the resulting solution, which includes the determination of the optimal switching time from one arc to the next one, violates another constraint, then the last two arcs are pieced together with the arc corresponding to the new violated constraint, and we re-solve the problem with the three arcs pieced together. The process is repeated until the solution does not violate any other constraints.

From \eqref{eq:min}, the state equations \eqref{eq:state} and the constraints \eqref{eq:safety},\eqref{eq:lateral_constraint}, for
each vehicle $i\in\mathcal{N}_z(t)$ the Hamiltonian function with the state and control adjoined is
\begin{gather}
H_{i}\big(t, x(t),u(t)\big) = \frac{1}{2} u^{2}_{i} + \lambda^{p}_{i} \cdot
v_{i} + \lambda^{v}_{i} \cdot u_{i} %
\nonumber\\
+ \mu^{a}_{i} \cdot(u_{i} - u_{max})+ \mu^{b}_{i} \cdot(u_{min} - u_{i}) + \mu^{c}_{i} \cdot(v_{i} - v_{max}) \nonumber\\
+ \mu^{d}_{i} \cdot(v_{min} - v_{i}), \quad \forall i \in\mathcal{N}_z(t),\label{eq:16b}%
\end{gather}
where $\lambda^{p}_{i}$ and $\lambda^{v}_{i}$ are the co-state components, and
$\mu_i^{a},\mu_i^{b},\mu_i^{c}$ and $\mu_i^{d}$ are the Lagrange multipliers.

If the inequality state and control constraints \eqref{eq:constraints} are not active, we have $\mu^{a}%
_{i} = \mu^{b}_{i}= \mu^{c}_{i}=\mu^{d}_{i}=0$. Applying the necessary condition, the optimal control can be given 
\begin{equation}
u_i(t) + \lambda^{v}_{i}= 0, \quad i \in\mathcal{N}_z(t). \label{eq:17}
\end{equation}
From the Euler-Lagrange equations, we have $\lambda^{p}_{i}(t) = a_{i}$, and $\lambda^{v}_{i}(t) = -\big(a_{i}\cdot t + b_{i}\big)$. 
The coefficients $a_{i}$ and $b_{i}$ are constants of integration corresponding to each vehicle $i$. From \eqref{eq:17} the optimal control input (acceleration/deceleration) as a function of time, and the corresponding state trajectories are given by
\begin{gather}
u^{*}_{i}(t) = a_{i} \cdot t + b_{i}, ~ \forall t \ge t^{0,z}_{i}. \label{eq:20}\\
v^{*}_{i}(t) = \frac{1}{2} a_{i} \cdot t^2 + b_{i} \cdot t + c_{i}, ~ \forall t \ge t^{0,z}_{i}\label{eq:21}\\
p^{*}_{i}(t) = \frac{1}{6}  a_{i} \cdot t^3 +\frac{1}{2} b_{i} \cdot t^2 + c_{i}\cdot t + d_{i}, ~ \forall t \ge t^{0,z}_{i}, \label{eq:22}%
\end{gather}
where $c_{i}$ and $d_{i}$ are constants of integration. The constants of integration $a_i$, $b_{i}$, $c_{i}$, and $d_{i}$ can be computed once at time $t^{0,z}_{i}$ using the initial and final conditions, and the values of the one of terminal transversality condition, i.e., $\lambda^{v}_{i}(t_i^{m,z})=0$. Similar results are obtained for the remaining cases. Due to space limitations, this analysis is omitted but may be found in \cite{Malikopoulos2017}. Note that, although the general framework of the upper-level time optimal CAV coordination and low-level energy optimal speed optimization is applicable to all conflict scenarios that we considered here, we incorporate the ``salient" features of the specific conflict scenarios by considering the control zone dimension $L_z$ and the merging zone dimension $S_z$ specific to the individual conflict scenario $z$ in Eqs. \eqref{eq:7} and \eqref{eq:8}. We also consider conflict scenario specific terminal conditions (e.g., $v_i(t_i^{m,z})$) to be applied in Eqs. \eqref{eq:20}-\eqref{eq:22} to derive the constants of integration of the control policy.
\section{Sequential Vehicle Testing}
\subsection{Simulation Environment}
To implement the control framework presented in the previous section, we use the microscopic multi-modal commercial traffic simulation software PTV VISSIM. We create a simulation setup replicating Mcity, and define a network consisting of three sources of bottleneck (one on-ramp merging, one speed reduction zone, and one roundabout), as shown in Fig. \ref{fig:mcity_corridor}. The corridor through which the Audi A3 travels has a length of 1500 $m$ length within Mcity. We select the CZ of each conflict scenario in such a way that no sharp turn or bends falls inside the controlled region. 
The vehicle dynamics parameters are calibrated by the data collected through driving the vehicle inside the Mcity facility. We consider different level of traffic congestion for investigating the robustness of the controller. A detailed analysis can be found in \cite{Zhao2019CCTAa}. In this paper, we consider the medium to high congestion case for simulation. We employ the desired traffic congestion by modifying the vehicle flow per hour per lane in the congestion route as illustrated by the black trajectories in Fig. \ref{fig:mcity_corridor}. Table \ref{tab:sim_param} contains the essential parameter required for setting up the simulation environment.


\begin{table}[!htb]
\fontsize{8}{10}\selectfont
\centering
\caption{VISSIM simulation parameters.}\label{tab:sim_param}
\begin{tabular}{| L{0.75\columnwidth-2\tabcolsep-1.2\arrayrulewidth} | L{0.25\columnwidth-2\tabcolsep-1.2\arrayrulewidth} |  }
\hline
\textbf{Vehicle Parameters} &                               \\ \hline
{Maximum Acceleration [$m/s^2$]} & 1.5                 \\  \hline
{Maximum Deceleration [$m/s^2$] }    & 3.0                   \\  \hline
{Safe Time Headway [$s$] }    & 1.2                   \\  \hline
\textbf{Traffic Network} &                               \\ \hline
{Corridor Length [$m$] }    & 1500                   \\  \hline
{Control Zone Length [$m$]} & 100                  \\  \hline
{SRZ Length [$m$]} & 125                  \\  \hline
\textbf{Speed Limit}    &                    \\  \hline
{On-Ramp Merging [$mph$] }    & 40                   \\  \hline
{Speed Reduction Zone [$mph$] }    & 18.6                   \\  \hline
{Roundabout [$mph$] }    & 25                   \\  \hline
\textbf{Traffic Flow }    &                    \\  \hline
{Main Route [$vph/lane$] }    & 500                   \\  \hline
{Highway [$vph/lane$] }    & 800                   \\  \hline
{Speed Reduction Zone [$vph/lane$] }    & 1300                   \\  \hline
{Roundabout [$vph/lane$] }    & 700                   \\  \hline
	\end{tabular}
	\par
   \vspace{-0.15\skip\footins}
   \renewcommand{\footnoterule}{}
\end{table}

\begin{figure}[t]
\centering
\includegraphics[width=3.5in]{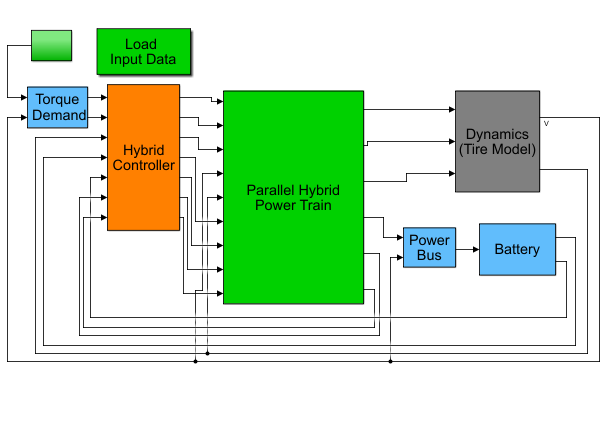} 
\caption{VESIM model for equivalent fuel consumption calculation (MPGe) of Audi A3 e-tron plugin hybrid electric vehicle.}%
\label{fig:vesim}%
\end{figure}

To evaluate the effectiveness of the proposed optimal vehicle dynamics control, we consider two different cases:

\begin{enumerate}

\item \textbf{Baseline scenario:} We construct the baseline scenario considering all the vehicles to be human-driven and without any V2V/V2I communication capability. The vehicles subscribe to the VISSIM built-in Wiedemann car following model \cite{wiedemann1974} to emulate the driving behavior of real human driven vehicles. We adopt priority based (yield/stop) traffic movement for the roundabout and on-ramp merging conflict scenarios. 

\item \textbf{Optimal VD controlled scenario:} In the optimal scenario, all the vehicles are CAVs, and communicate with each other inside the CZ through V2V/V2I communication. Therefore, they can plan their optimal path inside the CZ avoiding any lateral, or rear-end collision and optimize their individual travel time and fuel efficiency. We do not consider priority-based movement for the conflict scenarios as in the baseline case. We consider three isolated coordinators with a CZ of  $100\,m$ for each conflict zone (see Fig. \ref{fig:mcity_corridor}). For the uncontrolled paths in-between the CZ, the vehicles reverts back to the Wiedemann car following model \cite{wiedemann1974} to traverse their respective routes. To apply the proposed VD controller, we need to override the built-in car following module and associated attributes of the simulation environment. VISSIM provides various APIs as add-on modules to integrate VISSIM with user’s own applications. 
    
\end{enumerate}
Note that, the simulation in VISSIM only outputs the speed profile of the vehicle in the network, and does not yield any fuel consumption result specific to the Audi A3. Since, the Audi A3 is PHEV, fuel consumption models for conventional vehicles are not applicable. Therefore, we adopt a hybrid electric vehicle simulation model VESIM model, reported in \cite{Malikopoulos:2006aa} and references therein, calibrated appropriately to emulate the fuel consumption of the Audi A3 (see Fig. \ref{fig:vesim}). Due to the combined contribution of the internal combustion engine and the motor of the Audi A3, the VESIM model calculates the miles-per-gallon of gasoline equivalent (MPGe) according to the EPA standard.
By feeding the baseline and the VD controlled speed profiles to the VESIM model, we quantify the fuel consumption of the Audi A3, and evaluate the performance of the VD controller at different conflict scenarios.

\subsection{Hardware-in-the-Loop (HIL) Test}
We investigate the performance of the VD controller in a HIL environment through the chassis-dyno setup.  In our application, only a limited ECU variables of the Audi A3 can be monitored and, if necessary, bypassed. The rest of the variables are solely controlled by the Audi A3’s software system and cannot be accessed externally. To bypass the control input (acceleration) expected from driver pedal, we override the ECU variables responsible for the Audi A3's cruise controller. We also override the ECU variables relevant to the torque-regeneration module of the integrated motor-generator (IMG) unit to introduce braking (deceleration) force when required. That is, we use the brake regeneration functionality of the IMG unit of the Audi A3 to apply braking force when decelerating. However, the braking power generated by the IMG unit in this way is not as accurate as the typical hydraulic/mechanical brakes. Note that, applying brake through IMG is only possible when the state of the charge (SOC) of the battery is below a certain threshold to allow charging. If the battery SOC is above this certain threshold value, braking through IMG is not possible. The bypass of the aforementioned ECU variables has been realised by the ETAS ES910 rapid prototyping device. We calibrate the override coefficients in order to improve the tracking performance of the vehicle. Afterwards, we feed the speed profiles generated by simulation environment for both the baseline and optimal controlled case. We consider three different traffic volumes (e.g., low, medium, and high) to investigate the robustness of the VD controller performance, and record and analyze the fuel consumption data. 

\subsection{Bench Test of Integrated System}

\begin{figure}[h]
\centering
\includegraphics[width=3.5in]{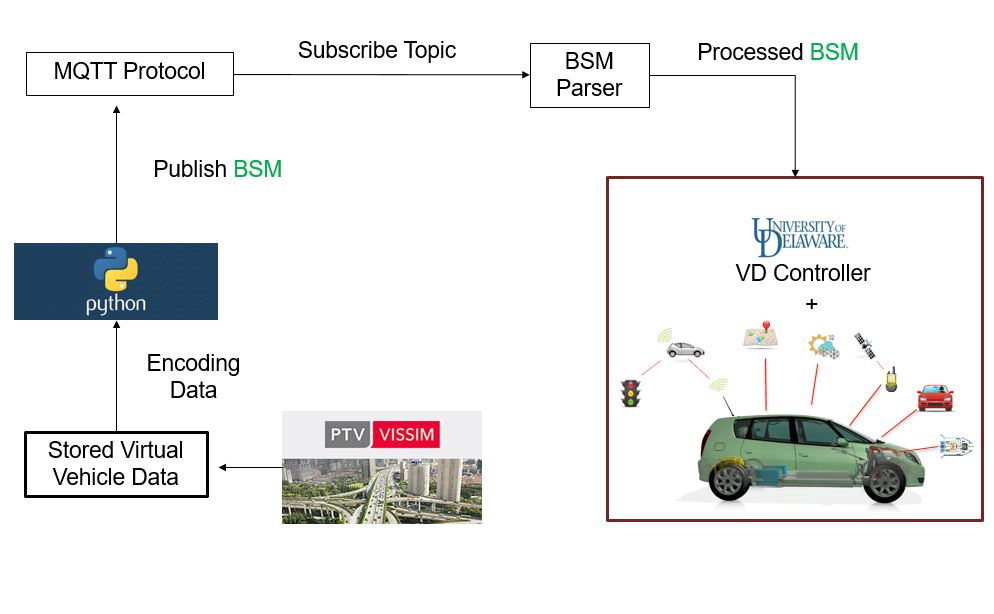} 
\caption{Bench test model emulating the virtual reality based field test.}%
\label{fig:benchTest}%
\end{figure}

We develop an IoT-based virtual reality framework to implement and debug the proposed VD controller on the bench. This particular bench test acts as a prerequisite to the actual field test, and performs general system validation of the integrated VD controller. In Fig. \ref{fig:benchTest}, we illustrate the proposed framework, which enables system validation before the field test. The framework we employ here can be carried out in a single computer without requiring any additional hardware or devices, which significantly reduces the setup and execution time.
The bench-test architecture in Fig. \ref{fig:benchTest} can be subdivided predominantly into two platforms, namely the simulation environment and the emulated head unit of the Audi A3. The connectivity of the platforms is established using the MQTT connectivity protocol. First, virtual vehicle information generated from the aforementioned VISSIM simulator is stored offline. To mimic the data transfer process via the V2X communication framework, a python script reads the stored virtual vehicle information, parses them into suitable Basic Safety Message (BSM) strings, and publishes them sequentially to a MQTT topic in the BSM format as expected during the actual field test. An MQTT subscription thread runs in the emulated head unit, which subscribes to the specific MQTT topic to receive the incoming BSMs of virtual vehicles. Once received, the head unit reads the published BSMs instantaneously, and uses its hard-coded decoder algorithm to parse and assign the necessary information to its internal memory. Afterwards, the head unit passes these stored information to the integrated VD controller to compute the optimal speed profile online. Note that, this is an open-loop process; namely, the test-vehicle information is not transmitted back to the simulation environment as the simulator works offline. With some simple modification, the process could be transformed into a closed loop one where the BSMs from the test vehicles is transmitted back to the simulator, and the virtual vehicles inside the simulator could react to the state of the test vehicle. 

\section{Virtual Reality Based Field Test}
\begin{figure}[htb]
\centering
\includegraphics[width = 3.5in]{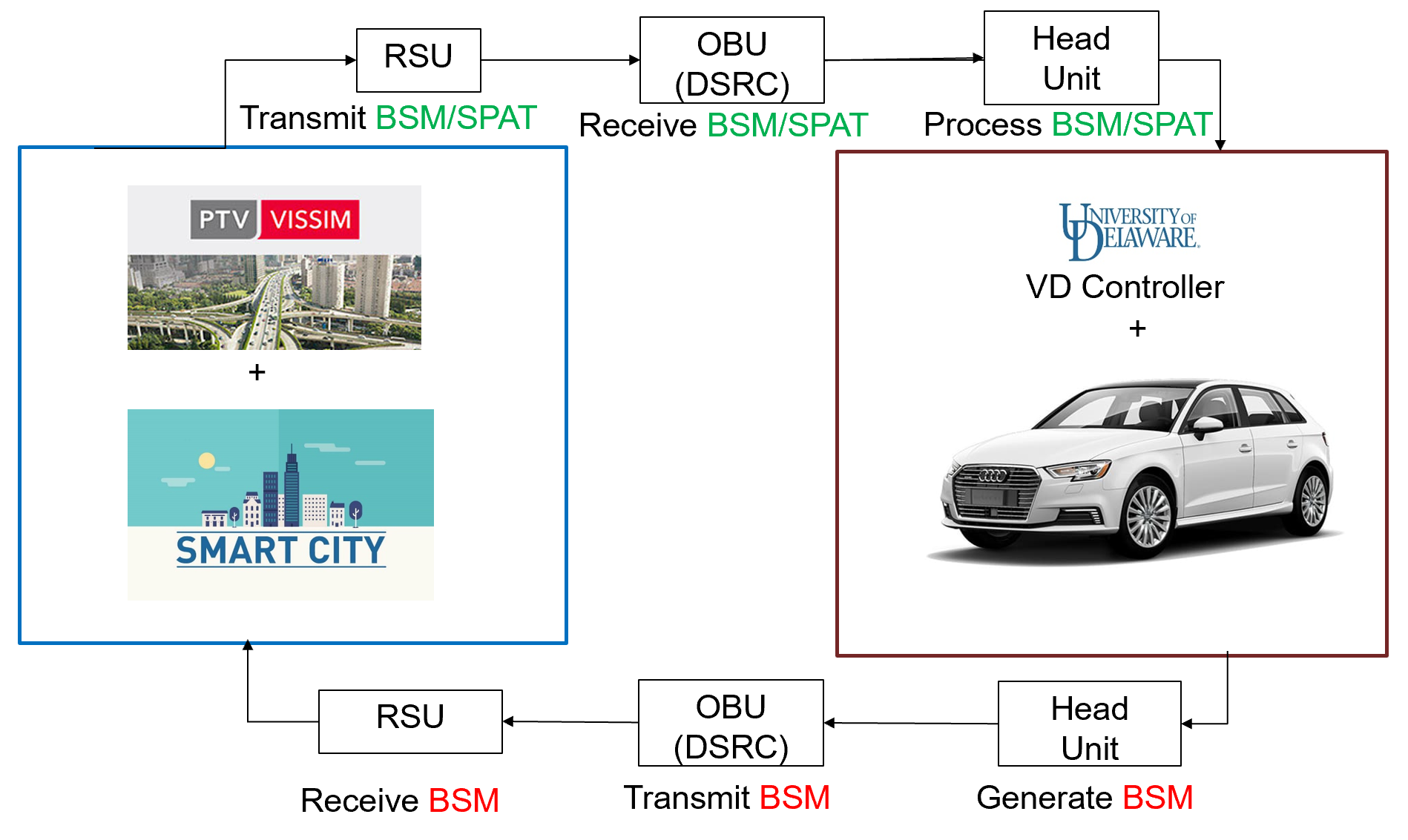} 
\caption{Workflow of virtual reality based V2X enabled field test in Mcity.}%
\label{fig:test_AR}%
\end{figure}

In this paper, we employ the general architecture of the virtual reality based test environment of Mcity, shown in Fig. \ref{fig:test_AR}. The test environment is composed of three main subsystems as follows:
\begin{enumerate}
    \item \textbf{Test Vehicle:} The Audi A3 is placed inside Mcity and is retrofitted with an on-board unit (OBU). The OBU considered here is a dedicated short range communication (DSRC) device from COHDA Wireless (model MK5). A GPS device (model BU-353-S4) is used to get the latitude and longitude of the Audi A3 during the test. The head unit of the Audi A3 is emulated by a Linux machine, which is connected to the Audi A3's ECU through
    ETAS ES-910 rapid prototyping device.
    \item \textbf{Infrastructure Equipment:} The infrastructure equipment includes roadside units (RSUs), traffic signal controllers and vehicle detectors. The RSUs act as two-way communication channel between the Audi A3 at the test facility, and the simulation platform located at the main control center.
    \item \textbf{Simulation Platform:} The virtual environment is provided by a simulation platform located at the main control center. The simulator can generate virtual vehicles and provide their information to the real-world structures. The simulation platform uses the commercial software PTV-VISSIM \cite{ptv2014ptv}, where different test scenarios can be saved as different projects. We construct the VESIM model based on the parameters presented in Table \ref{tab:sim_param}.
    The VISSIM  APIs, namely the SignalControl.DLL, DriverModel.DLL, and COM interfaces are used for interactions with the real-world environment and the simulation managing application. We modify the DriverModel.DLL to integrate the proposed VD controller applicable to all the virtual vehicles within the simulation platform. To enable the VISSIM simulator to receive and transmit BSMs, we embed the DriverModel.DLL with appropriate BSM encoder and decoder. At each time step, each virtual vehicle communicates with one another internally, and derives their own optimal speed. Information from simulated virtual vehicles are then encoded, and sent out by the DriverModel.DLL to the real-world Audi A3.
\end{enumerate}

The work flow of the vehicle testing procedure through virtual reality in Mcity is shown in Fig. \ref{fig:test_AR}. The Audi A3 communicates with the RSUs through the OBU, and receives basic safety messages (BSMs) coming from the simulation platform, and the signal phase and timing (SPaT) messages coming from the infrastructure. The Audi A3 computes its optimal desired speed, and passes this speed recommendation onto its ECU through CAN. At the same time, the Audi A3 broadcasts its state information (GPS location and speed) back to the RSUs through the OBU. Note that the OBU considered here is a DSRC device which transmits messages to/from the RSU. The RSU's data processor receives and processes the incoming data from the Audi A3, and sends the processed information to the simulation platform at the main control center. All of the traffic network and infrastructure attributes of the vehicle test facility are virtually designed and modeled in a simulation environment. Virtual vehicles are generated and their trajectories are updated based on the preset simulation setting and the information received from the Audi A3. Finally, the virtual vehicles in the simulation platform broadcast their information as BSMs to the RSUs, which is then received by the OBU of the Audi A3 and the process repeats itself. Through this communication framework, both the virtual vehicles in the simulation environment and the Audi A3 in Mcity can interact with each other real-time throughout the test session. Therefore, the behavior of both virtual and Audi A3 are completely synchronized. However, there are some drawbacks to this approach. First, significant efforts have to be taken for creating a proper simulation environment and communication framework between the real and virtual world. Second, preparing the Audi A3 for such an environment might not be straightforward, and would require significant involvement of the test facility itself. Lastly, the developed framework might not be easily transferable to other vehicle test facilities without significant modifications.


\subsection{Audi A3 System Architecture}
The system architecture inside the HU of the Audi A3 consists of the following subsystems,

\begin{enumerate}
    \item \textbf{GPS thread:} Collects and transmits the vehicle's latitude, longitude, and velocity information.
    \item \textbf{BSM thread:} Collects the incoming BSMs from the DSRC device, parses them into usable values, and transmits them to the main process.
    \item \textbf{SPaT thread:} Collects the incoming SPaT information from the DSRC device, parses them into usable information, and transmit them to the main process.
    \item \textbf{Mainframe:} Based on the connectivity protocol, the mainframe can have a single or multiple subscription threads to receive the incoming information from the previous threads. The mainframe houses the integrated VD controller which uses the incoming information. Once the VD controller outputs an optimal recommended speed, the mainframe transmits it to the vehicle through CAN.
\end{enumerate}

Based on the choice of connectivity framework, we develop two architectures as shown in Fig. \ref{fig:mqtt_protocol} and Fig. \ref{fig:socket_protocol}.

\subsubsection{MQTT and Socket Based Connectivity Protocols}
\begin{figure}
\vspace{1.5cm}
\centering
\includegraphics[width=3.5in]{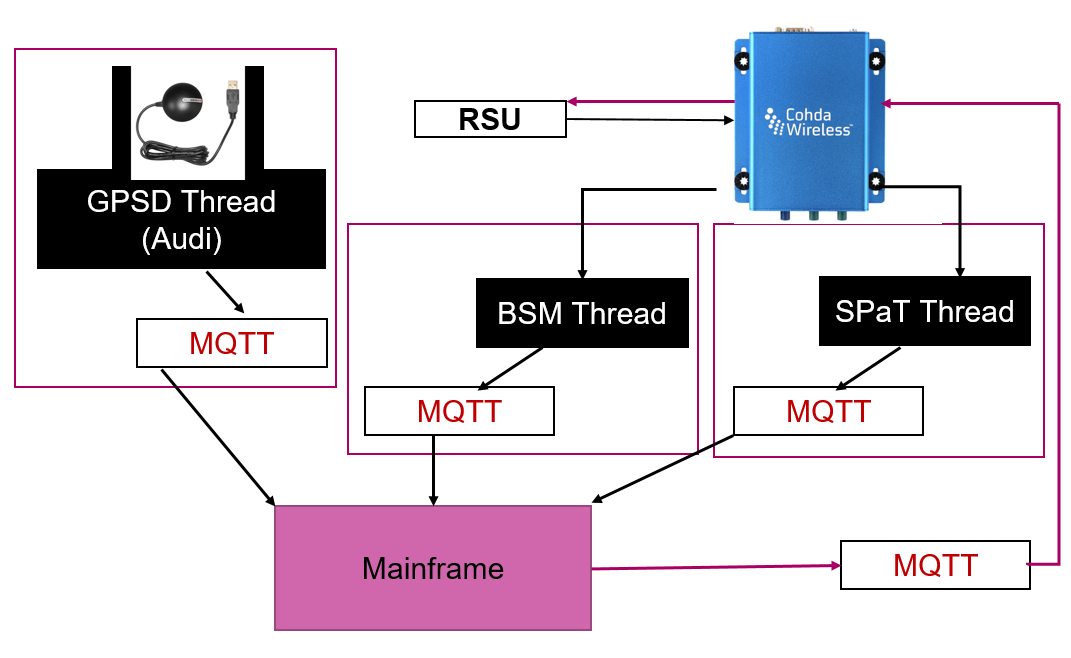} 
\caption{Workflow of MQTT Connectivity for individual processes.}%
\label{fig:mqtt_protocol}%
\end{figure}

MQTT is a M2M/IoT connectivity protocol used to publish/subscribe messaging transport. 
In the architecture illustrated in Fig. \ref{fig:mqtt_protocol}, the GPS, BSM, and SPaT information acquisition and transmission threads run as individual processes. Each of these processes are eventually connected to the mainframe of the head unit through basic MQTT protocol. The MQTT protocol is realized by an inner publish and subscription thread, which means that for every data transmission, the sender has to publish the data to a specific topic, which is, in turn, retrieved by the receiver through subscription to that specific topic. Due to requirement of several publish/subscription thread, the resulting system structure becomes complicated with additional data transmission and processing delays. The system delay can have a negative impact as it may prevent the real-time synchronization of virtual vehicles in the simulation platform and the Audi A3.

\begin{figure}[h]
\vspace{0.5cm}
\centering
\includegraphics[width=3.5in]{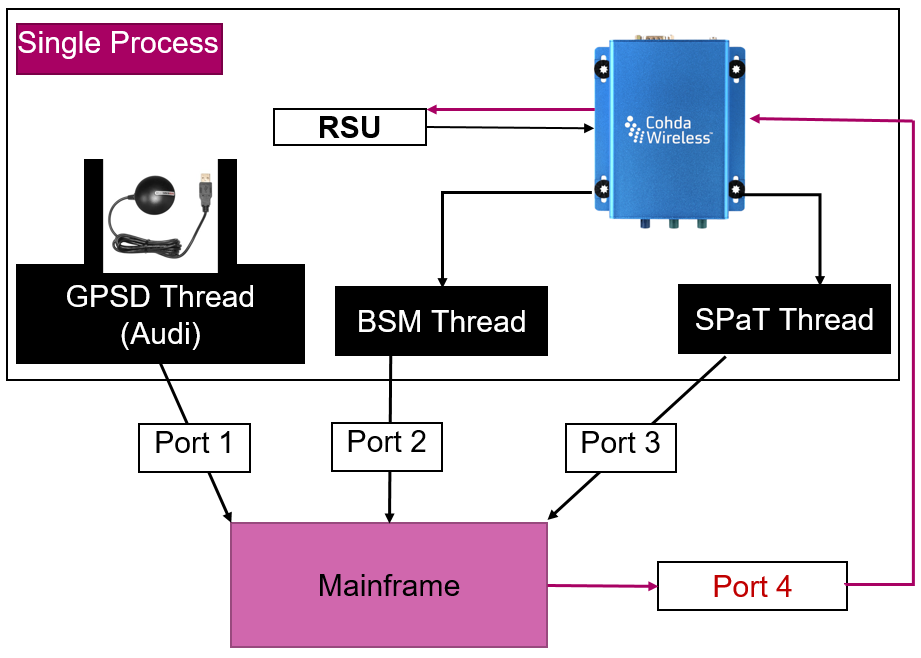} 

\caption{Workflow of socket based connectivity for the unified architecture.}%

\label{fig:socket_protocol}%
\end{figure}

To reduce the system delay associated with the MQTT based architecture, we adopt an improved version of the communication framework (Fig. \ref{fig:socket_protocol}). In this updated version, we  unify all the different threads (GPS, BSM and SPaT) under a single process. We completely eliminate the publish/subscription-based MQTT protocol, and use significantly faster connection through direct sockets. A socket is essentially defined as a connection having a unique IP address and relevant port number. By tapping into the IP address of the DSRC device, and using different ports to handle the BSM and SPaT information, we have greatly simplified the data acquisition and transmission procedure. The resulting framework has been tested to be significantly faster than the older version. For example, the head unit's data transmission frequency with the MQTT protocol is 1 Hz, whereas it increased to greater than 100 Hz with the socket-based protocol. The aforementioned modifications have been implemented to collect the optimal control field test data for the on-ramp merging, SRZ and the roundabout presented in this paper.

\subsection{VD Controller Integration}

\begin{figure}[h]
\vspace{0.5cm}
\centering
\includegraphics[width = 3.5in]{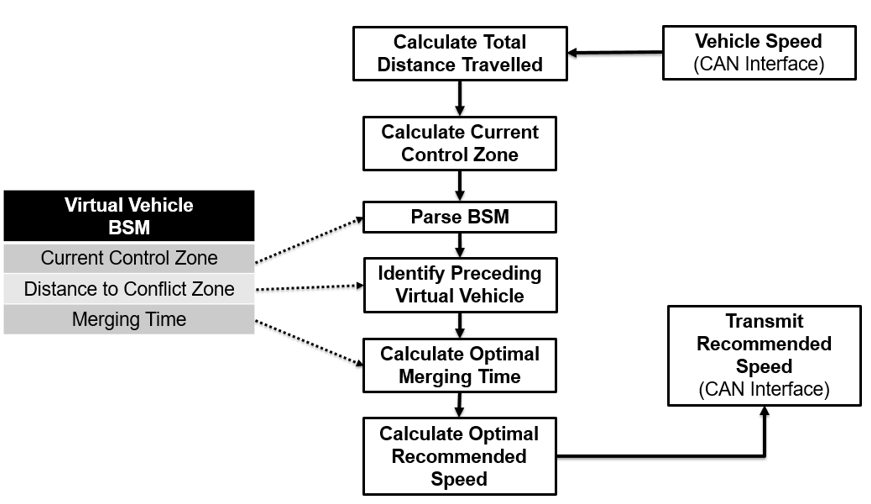} 
\caption{Workflow of VD Controller implementation within the head unit of Audi A3.}%
\label{fig:vd_controller}%
\end{figure}

We integrate the VD controller algorithm in the emulated head unit of Audi A3, and link it up with the BSM and SPAT encoder/decoder. Therefore, the VD controller can receive the essential information necessary for Audi A3's optimal speed calculation. To calculate the optimal speed trajectory, the Audi A3's head unit has to find the answers to the following sequential questions:
\begin{enumerate}
    \item Which control scenario does the test vehicle currently belong to? 
    \item Which particular virtual vehicle is acting as the putative leader at this specific control scenario?
    \item What is the merging time of that particular putative leader?
\end{enumerate}
Note that, the answer to the above questions can not be obtained directly from the incoming virtual vehicle information available from the regular BSM string as specified in the SAE-J2735 standard \cite{sae_J2735}. Therefore, we override some of the data sent by the regular BSM string to include the specific data useful to answer the aforementioned questions. For example, we override the regular data to be transmitted through the BSM's elevation, vehicle length and vehicle width variables to contain respectively the vehicle's merging time, vehicle's distance to the conflict zone, and vehicle's current CZ information. These updated BSM information is fed to the head unit of the Audi A3 where the optimal control input is computed through a workflow illustrated in Fig. \ref{fig:vd_controller}.

From the initialization of the test session, the Audi A3 continuously computes its total distance traveled by integrating its current speed over each time step. Based on the total distance traveled and the hard-coded location values of the conflict zones, the Audi A3 determines whether it has entered the CZ of a particular conflict scenario. Once the Audi A3 is within a specific CZ, it uses the CZ information from the incoming BSM strings of the virtual vehicles and sorts out the virtual vehicles pertaining only to that particular CZ. Afterwards, the Audi A3 sorts the virtual vehicles pertaining to the corresponding CZ by the distance to conflict zone information to find the preceding virtual vehicle. Finally, the Audi A3 uses the merging time, $t_{i-1}^{m,z}$ of the determined preceding virtual vehicle to calculate its own optimal speed input using Eq. \eqref{eq:21}. If there are no preceding virtual vehicles at any time $t$, the Audi A3 computes the optimal speed trajectory based on its estimated arrival time. The optimal speed input is transmitted to the vehicle's ECU via CAN-Bus. 

%
%
\section{Result and Discussion}
\subsection{Simulation Result}

\begin{figure}[h]
\centering
\includegraphics[width=3.5in]{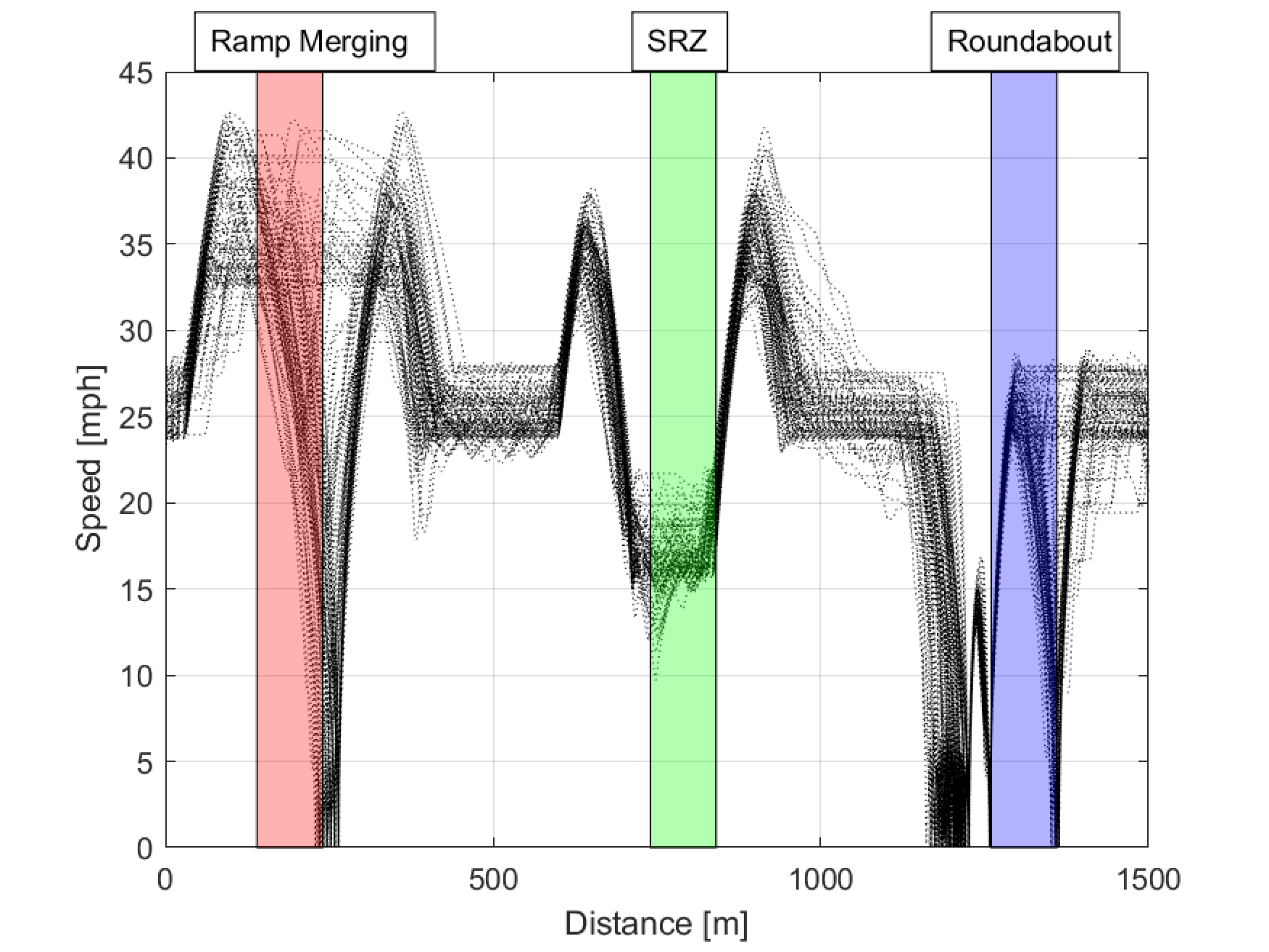}
\caption{Speed profile of baseline vehicles in simulation for high traffic volume.}%
\label{fig:simBaseline}%
\end{figure}

\begin{figure}[h]
\centering
\includegraphics[width=3.5in]{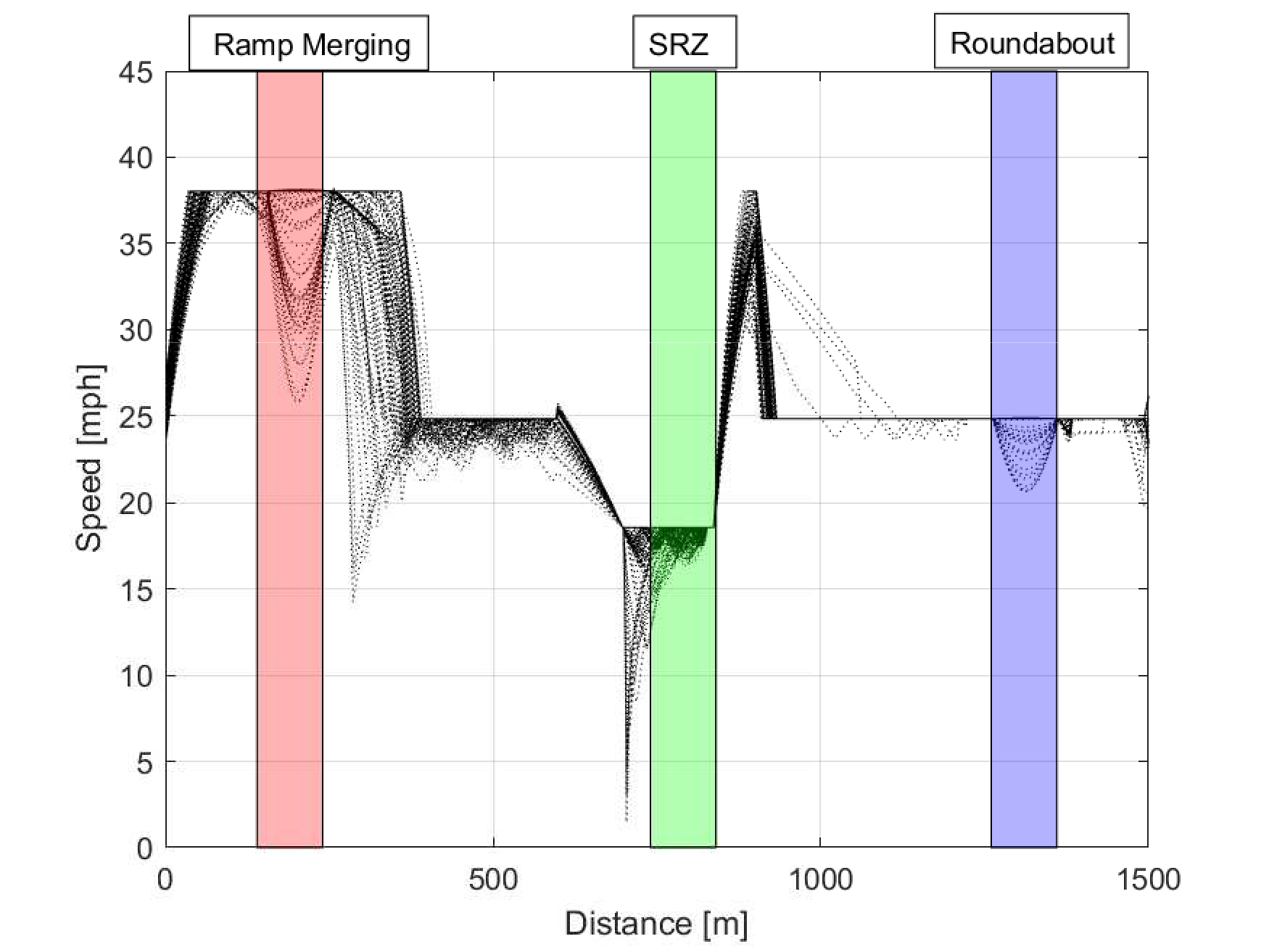} 
\caption{Speed profile of optimal controlled CAVs in simulation for high traffic volume.}%
\label{fig:simOptimal}%
\end{figure}

\begin{enumerate}

\item \textit{Baseline Scenario:} The speed profiles of human-driven vehicles (baseline scenario) with Wiedemann human driver model \cite{wiedemann1974} in the test route are shown in Fig. \ref{fig:simBaseline}. 
We notice significant stop-and-go driving at the on-ramp merging and roundabout conflict zones because the conventional human-driven vehicles have to yield to the incoming main-road vehicles. The human-driven vehicles can either show stop-and-go behavior, if the main road is very congested, or have a free merging, if the main road is empty. Moreover, we observe an increase in vehicle speed at the region before the SRZ. The SRZ is located near the end of a straight segment of Mcity (see Fig. \ref{fig:mcity_corridor}) which allows the vehicles to gain speed. Since the human-driven vehicles do not have any information regarding the upcoming SRZ or the slower moving vehicles in the SRZ, they pick up speed until the entry of the SRZ and starts decelerating at the entry of the SRZ.
We also observe that the stop-and-go driving behavior at each conflict scenario has negative implications in upstream and downstream of the route. Due to the presence of high speed variation inside the CZ, backward propagating traffic waves may result, which impacts the energy efficiency of the vehicles travelling outside the CZ. 

 \item  \textit{Optimal Controlled Scenario:} In Fig. \ref{fig:simOptimal}, we show the speed profiles of simulated CAVs under the optimal VD controller. We observe that the CAVs can be coordinated using the VD controller to space themselves out in such a way that they can pass through the conflict zones without stop-and-go driving. The VD controller smooths the traffic flow, and eliminates traffic congestion. We also note that in the SRZ, the speed of the CAVs is harmonized contributing to some additional benefits in fuel consumption. In this case, the CAVs know beforehand the state of the previous CAVs approaching the SRZ. Therefore, the CAVs can adjust their speed inside the control zone of the SRZ in such a way that they have smooth entry at the SRZ and negate the backward propagating traffic wave. Moreover, we observe that the driving behavior of the CAVs located upstream and downstream of each CZ are streamlined, which adds additional energy efficiency benefit to the vehicles travelling outside the CZ.
\end{enumerate}

Table \ref{tab:simResult} summarizes the performance evaluation of the VD controller in terms of fuel consumption, and shows the average improvements in each conflict zone in terms MPGe of the CAVs under medium and heavy traffic volumes. The human-driven vehicles in the baseline simulation scenario exhibit stop-and-go driving, which in turn increases the transient engine operation \cite{Malikopoulos2008a}. On the other hand, the optimal speed profiles of the VD controller in Fig. \ref{fig:simOptimal} reduce the CAVs' transient engine operation by eliminating the stop-and-go driving behavior. 


\begin{table}[!htb]
\fontsize{8}{10}\selectfont
\centering
\caption{Summary of fuel consumption [MPGe] comparison in simulation.}\label{tab:simResult}
\begin{tabular}{| L{0.325\columnwidth-2\tabcolsep-1.2\arrayrulewidth} | L{0.225\columnwidth-2\tabcolsep-1.2\arrayrulewidth} | L{0.225\columnwidth-2\tabcolsep-1.2\arrayrulewidth} | L{0.225\columnwidth-2\tabcolsep-1.2\arrayrulewidth} | }
\hline
           \textbf{Conflict Scenario} & \textbf{On-Ramp}         & \textbf{SRZ}         & \textbf{Roundabout}                      \\ \hline
\textbf{Improvement (MPGe)[\%]} & 15.6       & 21.2      & 35.3             \\  \hline
	\end{tabular}
	\par
   \vspace{-0.15\skip\footins}
   \renewcommand{\footnoterule}{}
\end{table}

\begin{table}[!htb]
\fontsize{8}{10}\selectfont
\centering
\caption{Summary of travel time comparison in simulation.}\label{tab:simResult_time}
\begin{tabular}{| L{0.325\columnwidth-2\tabcolsep-1.2\arrayrulewidth} | L{0.125\columnwidth-2\tabcolsep-1.2\arrayrulewidth} | L{0.125\columnwidth-2\tabcolsep-1.2\arrayrulewidth} | L{0.225\columnwidth-2\tabcolsep-1.2\arrayrulewidth} | 
L{0.2\columnwidth-2\tabcolsep-1.2\arrayrulewidth} | }
\hline
           \textbf{Conflict Scenario} & \textbf{On-Ramp}         & \textbf{SRZ}         & \textbf{Roundabout}  & \textbf{Corridor}                       \\ \hline
\textbf{Avg. Travel Time (Baseline) [s]} & 24.8       & 19.82   & 18.71  & 176.9             \\  \hline
\textbf{Avg. Travel Time (VD Controlled) [s]} & 15.9       & 19.78    & 13.4 & 131.1             \\  \hline
\textbf{Improvement [\%]} & 35.5       & 0.17    & 28.1 & 25.9             \\  \hline
	\end{tabular}
	\par
   \vspace{-0.15\skip\footins}
   \renewcommand{\footnoterule}{}
\end{table}

Table \ref{tab:simResult_time} reports on the average travel time of the CAVs to cross each of the conflict zones for both the baseline and the optimal VD controlled scenarios. We observe that, in the scenario corresponding to the VD controller, we have on average ~26$\%$ reduction in travel time for the whole corridor.

\subsection{HIL Test Result}
\begin{figure}[h]
\centering
\includegraphics[width=3.5in]{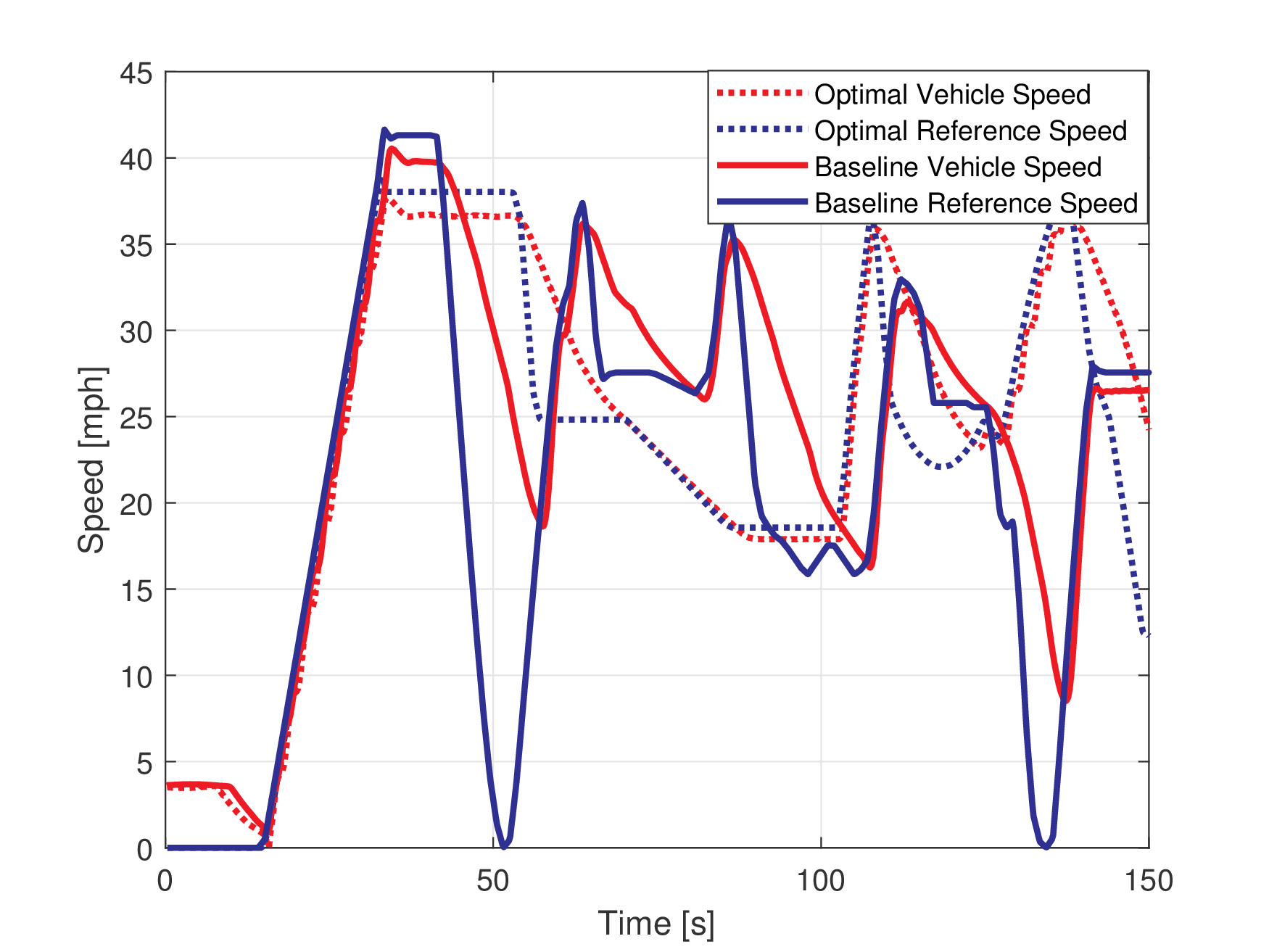} 
\caption{Reference and actual speed profile for baseline and VD controlled case under medium traffic volume.}%
\label{fig:HIL_mid}%
\end{figure}

\begin{figure}[h]
\centering
\includegraphics[width=3.5in]{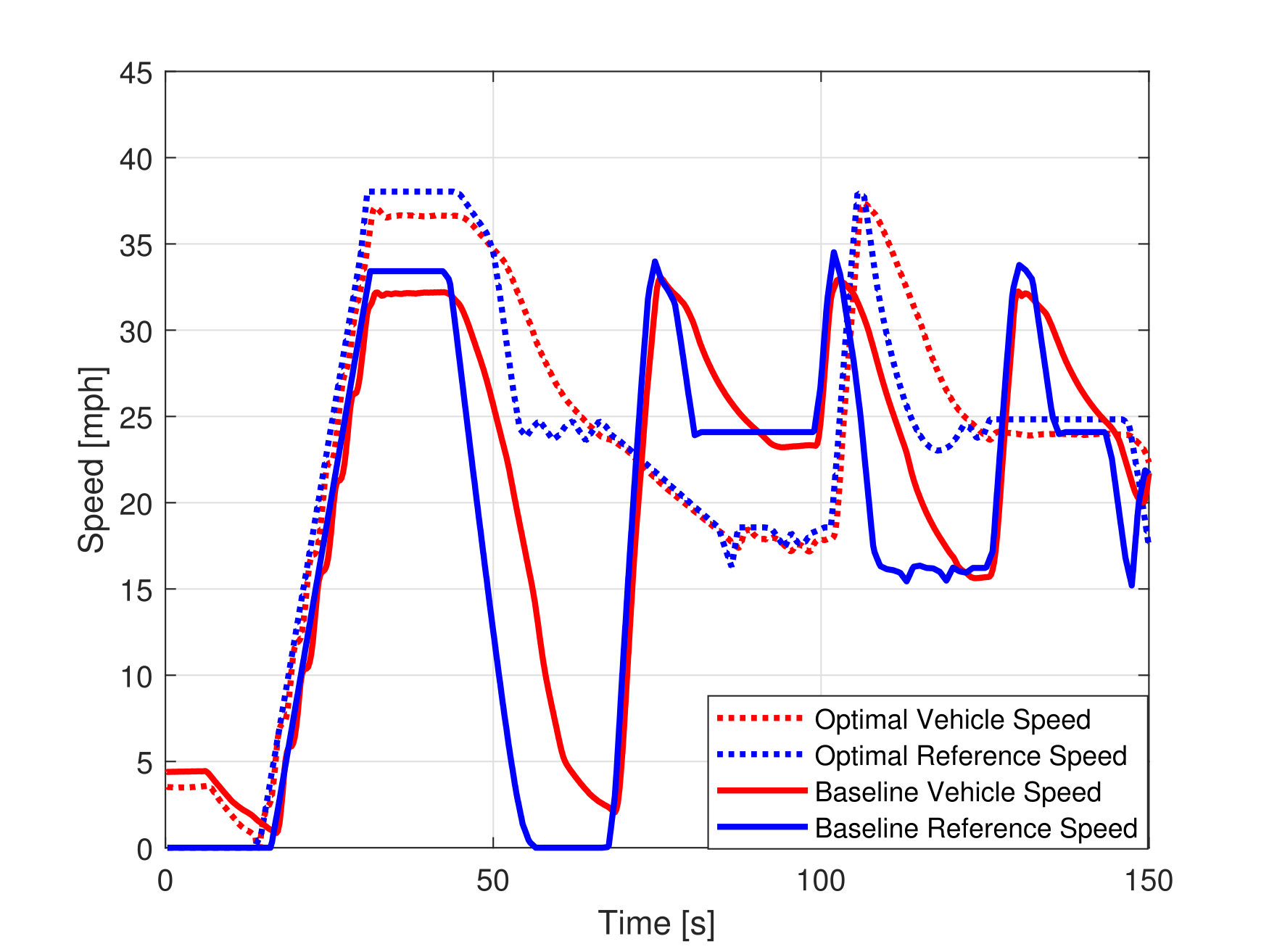} 
\caption{Reference and actual speed profile for baseline and VD controlled case under high traffic volume.}%
\label{fig:HIL_high}%
\end{figure}

The speed profiles generated by the hardware in the loop testing in the chassis-dyno setup are illustrated in Figs. \ref{fig:HIL_mid} and \ref{fig:HIL_high} for medium and high traffic volume respectively. We observe that the vehicle is able to trace the reference speed well while accelerating. While decelerating, we find discrepancy between the reference speed and the actual vehicle speed. As we mentioned earlier, we use the brake regeneration functionality of the IMG unit of the test vehicle to apply braking force when decelerating. However, the braking power generated by the IMG unit in this way is not as accurate as the typical hydraulic brakes.

Table \ref{tab:hilResult} summarizes the fuel consumption result for low, medium and heavy traffic volumes. 

\begin{table}[!htb]
\fontsize{8}{10}\selectfont
\centering
\caption{Summary of energy improvement in HIL Test.}\label{tab:hilResult}
\begin{tabular}{| L{0.325\columnwidth-2\tabcolsep-1.2\arrayrulewidth} | L{0.225\columnwidth-2\tabcolsep-1.2\arrayrulewidth} | L{0.225\columnwidth-2\tabcolsep-1.2\arrayrulewidth} | L{0.225\columnwidth-2\tabcolsep-1.2\arrayrulewidth} | }
\hline
           \textbf{Traffic} & \textbf{Baseline [MPGe]}         & \textbf{Optimal [MPGe]}         & \textbf{Improvement[\%]}                      \\ \hline
\textbf{Low }    & 27.7      & 37.1       & 33.7                \\  \hline
\textbf{Medium} & 36.6       & 49.0      & 34.2             \\  \hline
\textbf{High} & 36.1       & 50.8      & 41.0             \\  \hline
	\end{tabular}
	\par
   \vspace{-0.15\skip\footins}
   \renewcommand{\footnoterule}{}
\end{table}
\subsection{Field Test Result in Mcity}

\begin{figure}[h]
\centering
\includegraphics[width=3.5in]{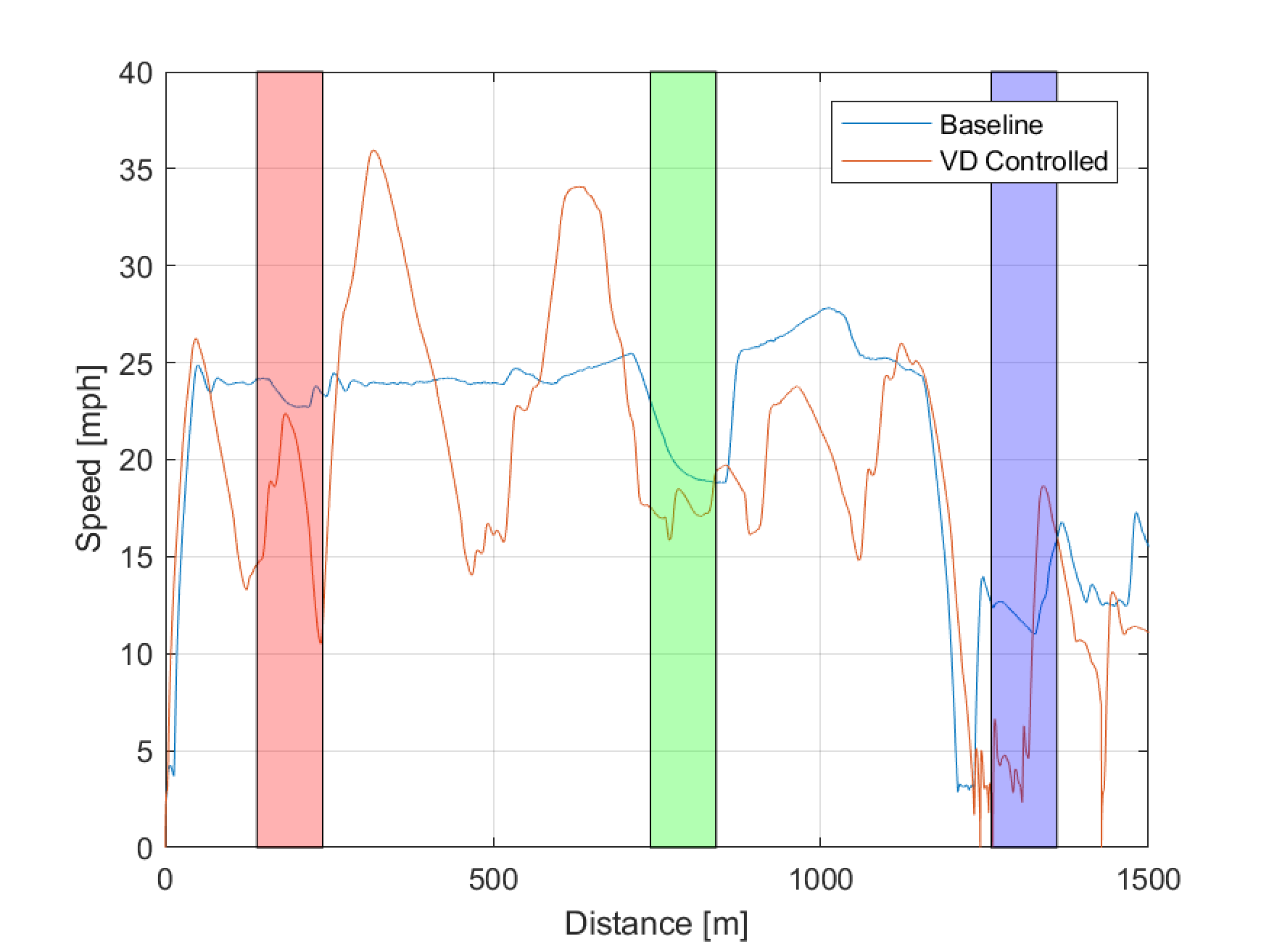} 
\caption{Speed profile of baseline vs optimal VD controlled case at the field test in Mcity.}%
\label{fig:fieldTest}%
\end{figure}

We summarize the field test results in Mcity in Table \ref{tab:resultFieldTest} and Table \ref{tab:resultfieldTest_time} which quantifies the Audi A3's performance in terms of fuel efficiency and travel time. The VD controlled Audi A3 shows improved fuel efficiency over the baseline scenario (human-driven Audi A3) exhibiting an overall 20 $\%$ improvement. A representative baseline and VD controlled speed profile of the Audi A3 from the field test in Mcity is shown in Fig. \ref{fig:fieldTest}. We observe that the optimal VD controlled speed profile is comparatively smoother than the baseline one.

The optimal VD controller also shows improvement in terms of travel time through each of the conflict zones. Essentially, by optimizing the speed profile, and eliminating the stop-and-go behavior at the bottlenecks, the VD controller achieves 17$\%$ reduction in travel time in the whole corridor. 

\begin{table}[!htb]
\fontsize{8}{10}\selectfont
\centering
\caption{Summary of MPGe Improvement in Field Test.}\label{tab:resultFieldTest}
\begin{tabular}{| L{0.325\columnwidth-2\tabcolsep-1.2\arrayrulewidth} | L{0.225\columnwidth-2\tabcolsep-1.2\arrayrulewidth} | L{0.225\columnwidth-2\tabcolsep-1.2\arrayrulewidth} | L{0.225\columnwidth-2\tabcolsep-1.2\arrayrulewidth} | }
\hline
            & \textbf{Optimal [MPGe]}         & \textbf{Baseline [MPGe]}         & \textbf{Improvement[\%]}                      \\ \hline
\textbf{On-Ramp}    & 57.3      & 47.5       & 20.6                \\  \hline
\textbf{SRZ} & 63.3       & 48.6      & 30.2             \\  \hline
\textbf{Roundabout} & 65.4       & 57.7      & 13.3             \\  \hline
\textbf{Total} & 64.8       & 54.0      & 20.0             \\  \hline
	\end{tabular}
	\par
   \vspace{-0.15\skip\footins}
   \renewcommand{\footnoterule}{}
\end{table}

\begin{table}[!htb]
\fontsize{8}{10}\selectfont
\centering
\caption{Summary of travel time comparison of field test.}\label{tab:resultfieldTest_time}
\begin{tabular}{| L{0.325\columnwidth-2\tabcolsep-1.2\arrayrulewidth} | L{0.125\columnwidth-2\tabcolsep-1.2\arrayrulewidth} | L{0.125\columnwidth-2\tabcolsep-1.2\arrayrulewidth} | L{0.225\columnwidth-2\tabcolsep-1.2\arrayrulewidth} | 
L{0.2\columnwidth-2\tabcolsep-1.2\arrayrulewidth} | }
\hline
           \textbf{Conflict Scenario} & \textbf{On-Ramp}         & \textbf{SRZ}         & \textbf{Roundabout}  & \textbf{Corridor}                       \\ \hline
\textbf{Avg. Travel Time (Baseline) [s]} & 34.3       & 19.4   & 15.3  & 218.2             \\  \hline
\textbf{Avg. Travel Time (VD Controlled) [s]} & 32.4       & 19.2    & 9.38 & 180.2             \\  \hline
\textbf{Improvement [\%]} & 5.6      & 1.0    & 38.8 & 17.4             \\  \hline
	\end{tabular}
	\par
   \vspace{-0.15\skip\footins}
   \renewcommand{\footnoterule}{}
\end{table}

\section{Conclusion}
\label{sec:concl}
In this paper, we presented a framework to validate a VD controller. We investigated the effectiveness of the controller systematically through a sequence of approaches, leading up to the virtual reality based real-world field test. We quantified the Audi A3's performance for baseline and optimal VD controlled scenarios, and conducted a comparative analysis. We concluded that the Audi A3 retrofitted with the optimal VD controller shows significant improvements in terms of fuel consumption and travel time. To validate the route-agnostic improvement of the proposed VD controller, we presented the performance metrics for individual conflict scenarios.

\section{Recommendation}
\label{sec:recom}
Ongoing research is addressing several issues faced during the field test procedure (e.g., road grade, state/control constraints, data latency etc.). 
Future efforts will address the VD controller's extension to the mixed traffic conditions (partial penetration of CAVs), where the human-driven vehicles and the CAVs share the network simultaneously at different composition.

\bibliographystyle{ieeetr} 
\bibliography{field_test_VD_ref}

\section{Contact Information}
A M Ishtiaque Mahbub, \newline
mahbub@udel.edu \newline

Andreas A. Malikopoulos, \newline
andreas@udel.edu 
\section{Acknowledgments}
This research was supported by ARPAE's NEXTCAR program under the award number DE-AR0000796.

\section{Definitions, Acronyms, Abbreviations}

\begin{table}[h]
\centering
\begin{tabular}{L{0.1\textwidth} L{0.33\textwidth}}
\textbf{BSM} & Basic Safety Message \\
\textbf{CAV} & Connected Automated Vehicle \\
\textbf{CAN} & Control Area Network \\
\textbf{CZ}  & Control Zone \\
\textbf{DSRC} & Dedicated Short Range Communication \\
\textbf{ECU} & Electronic Control Unit \\
\textbf{HU} & Head Unit \\
\textbf{HIL} & Hardware-in-the-Loop \\
\textbf{IMG} & Integrated Motor Generator \\
\textbf{MQTT} & Message Queuing Telemetry Transport\\
\textbf{MPGe} & Miles-Per-Gallon of Gasoline Equivalent \\
\textbf{MZ} & Merging Zone\\
\textbf{OBU} & On Board Unit \\
\textbf{PHEV} & Plugin Hybrid Electric Vehicle \\
\textbf{RSU} & Road Side Unit \\
\textbf{SRZ} & Speed Reduction Zone \\
\textbf{SPaT} & Signal Phase and Timing \\
\textbf{SOC} & State-of-Charge \\
\textbf{TCP/IP} & Transmission control protocol/Internet protocol \\
\textbf{VD} & Vehicle dynamics \\
\textbf{V2X} & Vehicle to everything \\
\textbf{V2V} & Vehicle to vehicle \\
\textbf{V2I} & Vehicle to infrastructure \\
\textbf{VPH} & Vehicle per hour \\
\end{tabular}
\end{table}



\end{document}